\documentclass[10pt]{article}

\usepackage{amsmath,amssymb,amsfonts,latexsym,theorem}
\theoremstyle{plain}
\newtheorem{conjecture}{Conjecture}

\newtheorem{theorem}{Theorem}
\newtheorem{corollary}{Corollary}
\newtheorem{proposition}{Proposition}
{\theorembodyfont{\rmfamily}
\newtheorem{remark}{Remark}

}

\def\qed{\hfill $\square$}

\title{\textbf{Pattern Containment and Combinatorial Inequalities}}

\author{
Alexander I. Burstein\\
Department of Mathematics\\
Iowa State University\\
Ames, IA 50011-2064\\
\texttt{burstein@math.iastate.edu}
}

\begin{document}

\maketitle

\begin{abstract}
We use a probabilistic method to produce some combinatorial
inequalities by considering pattern containment in permutations
and words.
\end{abstract}

If $\sigma\in S_n$ and $\tau\in S_m$, we say that $\sigma$
\emph{contains} $\tau$, or $\tau$ \emph{occurs} in $\sigma$, if
$\sigma$ has a subsequence order-isomorphic to $\tau$. In this
situation, $\tau$ is called a \emph{pattern}. Similarly, if
$\sigma\in[k]^n$ is a string of $n$ letters over the alphabet
$[k]=\{1,\dots,k\}$, and $\tau\in[l]^m$ is a map from $[m]$ onto
$[l]$ (i.e. $\tau$ contains all letters from 1 to $l$), then we
say that $\sigma$ contains the pattern $\tau$ if $\sigma$ has a
subsequence order-isomorphic to $\tau$. An \emph{instance} (or
\emph{occurrence}) of $\tau$ in $\sigma$ is a choice of $m$
positions $1 \le i_1 < \ldots < i_m \le n$, such that the
subsequence $(\sigma(i_1),\ldots,\sigma(i_m))$ is order-isomorphic
to $\tau=(\tau(1),\ldots,\tau(m))$.

Most of the work on pattern containment concentrated on pattern
avoidance, that is on characterizing and counting permutations
that contain no occurrences of a given pattern or a set of
patterns. Less attention has been given to counting the number of
times a given pattern occurs in permutations of a given size, in
particular, packing patterns into permutations (but see
\cite{AAHHS,P}, for example), and, to our knowledge, packing
patterns into words (where repeated letters are allowed) has not
yet been considered.

Here we consider pattern containment and use a simple
probabilistic fact (the variance of a random variable is
nonnegative) to produce nontrivial combinatorial inequalities.

\section{Patterns in permutations}

In this section, we consider permutation patterns contained in
other permutations.

\begin{theorem} \label{perm}
Let $\tau$ be a permutation of $\{0,1,\ldots,m\}$ and define
\[
[i,j]_m=\binom{i+j}{i}\binom{2m-i-j}{m-i}.
\]
Then for any nonnegative integer $m$ and any $\tau$ as above,
\begin{equation} \label{eq-perm}
\sum_{i,j=0}^{m}{[i,j]_m [\tau(i),\tau(j)]_m} \ge
\binom{2m+1}{m}^2.
\end{equation}
\end{theorem}

\begin{remark} \label{perm-paths}
Notice that $[i,j]_m$ is the number of northeast integer lattice
paths from $(0,0)$ to $(m,m)$ through $(i,j)$. Hence the left-hand
side is the number of pairs $(P,Q)$ of northeast integer lattice
paths $P:(0,0)\to(i,j)\to(m,m)$ and
$Q:(0,0)\to(\tau(i),\tau(j))\to(m,m)$ over all $(i,j)\in[0,m]^2$.
\end{remark}

\begin{remark} \label{ekhad}
The numbers $[i,j]_m$, $0\le i,j\le m$ have been found to have
other interesting properties as well. For example, Amdeberhan and
Ekhad \cite{AE} showed that
\[
\det([i,j]_m)_{0\le i,j\le m}=\frac{(2m+1)!^{m+1}}{(2m+1)!!},
\]
where $a!!=0!\cdot 1!\cdot 2!\cdot\ldots\cdot a!$.
\end{remark}

It is a well-known result \cite{HLP} that for any nondecreasing
subsequence $a_1\ge\cdots\ge a_n$ and a permutation $\varphi\in
S_n$, the sum $\sum_{i=1}^{n}a_ia_{\varphi(i)}$ attains its
maximum when $\varphi=id_n=12\ldots n$ and its minimum when
$\varphi=n(n-1)\ldots 1$. Now if we arrange $(m+1)^2$ numbers
$[i,j]_m$ ($0\le i,j\le m$) in nondecreasing order, \emph{there is
no permutation $\tau$ of $\{0,1,\ldots,m\}$ which reverses that
order} (other than in the trivial case $m=0$). For example, even
when  $m=1$, we have
\[
[0,0]_2=[1,1]_2=2>1=[0,1]_2=[1,0]_2,
\]
reversing $(2,2,1,1)$ gives $(1,1,2,2)$ and the estimate of
\cite{HLP} gives us the lower bound of
$2\cdot1+2\cdot1+1\cdot2+1\cdot2=8$, our estimate yields the lower
bound of $\binom{3}{1}^2=9$, while the left hand side is actually
equal to 10 for both $\tau=01$ (the identity) and $\tau=10$ (which
transposes $[0,0]_2$ and $[1,1]_2$ as well as $[0,1]_2$ and
$[1,0]_2$ in the above ordering).

The estimate in Theorem \ref{perm} appears to be stronger than
that of \cite{HLP}. For example, the lower bounds for $m=2,3,4,5$
are $75,792,8660,98876$, respectively, according to \cite{HLP},
while our lower bounds are $100,1225,15876,213444$, respectively.
In fact, as the following proposition shows, the lower bound of
\cite{HLP} can never be achieved in our case for $m>0$.

\begin{proposition} \label{prop1}
Arrange $(m+1)^2$ numbers $\{[i,j]_m \mid 0\le i,j\le m\}$ into a
nondecreasing order $a_1\ge\ldots\ge a_{(m+1)^2}$. A permutation
$\tau$ of $\{0,1,\ldots,m\}$ induces an equivalence class of
permutations $\varphi_\tau$ on the $a_i$'s (equivalence relation
being a permutation of equal elements). Then for any $\tau$,
reversal of the identity
$(m+1)^2((m+1)^2-1)\ldots21\notin\varphi_\tau$.
\end{proposition}

\noindent\textbf{Proof.} Suppose that there is a permutation
$\tau$ which induces an order-reversing permutation of the
$a_i$'s. Note that $[0,m]_m=[m,0]_m=1$ for any $m$, hence,
$[\tau(0),\tau(m)]_m$ must have the greatest value among all
$[\tau(i),\tau(j)]_m$. Note that
\[
[i,j]_m=[j,i]_m=[m-i,m-j]_m=[m-j,m-i]_m
\]
for any $i$ and $j$, so assume that $i\le j$. Then it is a
straightforward exercise to prove that
\begin{align*}
[i,j]_m&>[i-1,j]_m \quad \text{for $i>0$, and}\\
[i,j]_m&>[i,j+1]_m \quad \text{for $j<m$,}
\end{align*}
so for $0<i<j<m$,
\[
[i,i]_m>[i,j]_m>[0,j]_m>[0,m]_m=1.
\]
Similarly, we can assume that $i\le\lfloor\frac{m}{2}\rfloor$
(since $[i,i]_m=[m-i,m-i]_m$), then it is just as easy to see that
for any $i>0$
\[
[i,i]_m<[i-1,i-1]_m<\cdots<[0,0]_m=\binom{2m}{m}.
\]
Thus, for any $0\le i,j\le m$,
\[
1=[0,m]_m=[m,0]_m \le [i,j]_m \le [0,0]_m=[m,m]_m=\binom{2m}{m},
\]
and one of the two inequalities becomes an equality if and only if
$i,j\in\{0,m\}$. Hence, for our permutation $\tau$, we must have
$\tau(0)=\tau(m)=0$ or $\tau(0)=\tau(m)=m$, neither of which is
possible when $m\neq 0$. The resulting contradiction implies our
proposition. \qed

\bigskip

Finally, before we begin with the proof of Theorem \ref{perm}, let
us note that a permutation of summands in (\ref{eq-perm}) yields
the following corollary.

\begin{corollary} \label{cor1}
For any two permutations $\tau_1,\tau_2$ of $\{0,1,\ldots,m\}$ and
any $m\in\mathbb{N}$,
\[
\sum_{i,j=0}^{m}{[\tau_1(i),\tau_1(j)]_m[\tau_2(i),\tau_2(j)]_m}\ge
\binom{2m+1}{m}^2.
\]
\end{corollary}

Another immediate corollary is a consequence of the fact that
\[
[i,j]_m=\binom{2m}{m}\frac{\binom{m}{i}\binom{m}{j}}{\binom{2m}{i+j}}=
\binom{2m}{m}\{i,j\}_m, \text{ where }
\{i,j\}_m:=\frac{\binom{m}{i}\binom{m}{j}}{\binom{2m}{i+j}}.
\]
\begin{corollary} \label{cor2}
For any $m\in\mathbb{N}$ and any permutation $\tau$ of
$\{0,1,\ldots,m\}$,
\[
\sum_{i,j=0}^{m}{\{i,j\}_m\{\tau(i),\tau(j)\}_m}\ge
\left(\frac{2m+1}{m+1}\right)^2=\left(2-\frac{1}{m+1}\right)^2.
\]
\end{corollary}

Note that Corollary \ref{cor2} no longer holds if we substitute
$4$ on the right side of this inequality.

\noindent\textbf{Proof of Theorem \ref{perm}.} Consider $S_n$ as a
sample space with uniform distribution. Let $\tau\in S_m$
(notation-wise, it is more convenient if, in the proof, $\tau$ is
a permutation of $\{1,2,\ldots,m\}$), and let $X_\tau$ be a random
variable such that $X_\tau(\sigma)$ is the number of occurrences
of pattern $\tau$ in given permutation $\sigma\in S_n$. We will
show that our inequality follows from the fact that
\[
V\!ar(X_\tau)=E(X_\tau^2)-E(X_\tau)^2\ge 0
\]

We start by finding $E(X_\tau)$. Pick an $m$-letter subset $S$ of
$[n]=\{1,2,\ldots,n\}$ in $\binom{n}{m}$ ways. There is a unique
permutation $\tau(S)$ of $S$ which is order-isomorphic to $\tau$.
There are $m!$ equally likely permutations in which the elements
of $S$ can occur in $\sigma$, but we need only 1 of them, namely,
$\tau(S)$. Hence, $\tau(S)$ either occurs once or does not occur
in a given permutation $\sigma$. Therefore, the probability that a
random $\sigma$ contains $\tau(S)$ as a subsequence is $1/m!$. Let
$Y_{\tau(S)}$ be a random variable such that $Y_{\tau(S)}(\sigma)$
is the number of occurrences of $\tau(S)$ in $\sigma$. Since
\[
P\!\left(Y_{\tau(S)}(\sigma)=1\right)=\frac{1}{m!} \quad
\text{and} \quad
P\!\left(Y_{\tau(S)}(\sigma)=0\right)=1-\frac{1}{m!},
\]
we have $E(Y_{\tau(S)})=1/m!$. But this is true for any
$S\subseteq[n]$ such that $|S|=m$, and we have
\[
X_\tau=\!\!\!\!\!\!\sum_{S\subseteq[n],\,|S|=m}{\!\!\!\!\!\!Y_{\tau(S)}}
\]
hence,
\[
E(X_\tau)=\!\!\!\!\!\!\sum_{S\subseteq[n],\,|S|=m}{\!\!\!\!\!\!E(Y_{\tau(S)})}=
\frac{1}{m!}\binom{n}{m}.
\]

Next, we look at $E(X_\tau^2)$. We have
\[
E(X_\tau^2)=E\left(\sum_{S\subseteq[n],\,|S|=m}{\!\!\!\!\!\!Y_{\tau(S)}}\right)=
\!\!\!\!\!\!\sum_{\substack{S_1,S_2\subseteq[n]\\|S_1|=|S_2|=m}}
{\!\!\!\!\!\!E\left(Y_{\tau(S_1)}Y_{\tau(S_2)}\right)}.
\]

Of course, $Y_{\tau(S_1)}Y_{\tau(S_2)}=1$ if and only if both
$\tau(S_1)$ and $\tau(S_2)$ are subsequences of $\sigma$,
otherwise, $Y_{\tau(S_1)}Y_{\tau(S_2)}=0$.

Let $S=S_1\cup S_2$,  and $|S_1\cap S_2|=\ell$, so $|S|=|S_1\cup
S_2|=2m-\ell$. We can pick a subset $S\subseteq[n]$ in
$\binom{n}{2m-\ell}$ ways. Note that any such $S$ is
order-isomorphic to $[2m-\ell]=\{1,2,...,2m-\ell\}$. Therefore,
the number of permutations $\rho(S)$ of $S$ such that
$\rho(S)=\tau(S_1)\cup\tau(S_2)$ for some $S_1,S_2\subseteq S$,
$S_1\cup S_2=S$, is the same for any $S$ of cardinality $2m-\ell$
and depends only on $m$ and $\ell$.

Therefore, $E(X_\tau^2)$ is a linear combination of
$\left\{\binom{n}{2m-\ell}\,\mid\,0\le\ell\le m\right\}$ with
coefficients which are rational functions of $m$ and $\ell$. The
degrees in $n$ of both $E(X_\tau^2)$ and $E(X_\tau)^2$ are $2m$,
and the coefficient of $n^{2m}$ in $E(X_\tau)^2$ is $1/(m!)^4$. On
the other hand, $S=S_1\cup S_2$, $|S|=2m$ and $|S_1|=|S_2|=m$
imply that $S_1\cap S_2=\emptyset$, so $Y_{\tau(S_1)}$ and
$Y_{\tau(S_2)}$ are independent, and hence
\[
P\left(Y_{\tau(S_1)}Y_{\tau(S_2)}=1\right)=P\left(Y_{\tau(S_1)}=1\right)
P\left(Y_{\tau(S_2)}=1\right)=\left(\frac{1}{m!}\right)^2.
\]
Since the number of ways to partition a set $S$ of size $2m$ into
two subsets of size $m$  is $\binom{2m}{m}$, the coefficient of
$\binom{n}{2m}$ in $E(X_\tau^2)$ is $\binom{2m}{m}/(m!)^2$. Hence,
the coefficient of $n^{2m}$ in $E(X_\tau^2)$ is
\[
[n^{2m}]E(X_\tau^2)=\frac{1}{(2m)!}\frac{1}{(m!)^2}\binom{2m}{m}=
\frac{1}{(m!)^4},
\]
where $[x^d]P(x)$ denotes the coefficient of $x^d$ in a given
polynomial $P(x)$. But then
$[n^{2m}]E(X_\tau^2)=[n^{2m}]E(X_\tau)^2$, so
$\deg_n(V\!ar(X_\tau))\le 2m-1$, and hence,
$[n^{2m-1}]V\!ar(X_\tau)\ge 0$.

We have
\begin{multline*}
[n^{2m-1}]E(X_\tau)^2=[n^{2m-1}]\left(\frac{1}{m!}\binom{n}{m}\right)^2=\\
=\frac{2}{(m!)^2}\cdot[n^{m}]\binom{n}{m}\cdot[n^{m-1}]\binom{n}{m}=\\
=\frac{2}{(m!)^2}\cdot\frac{1}{m!}\cdot\left(-\frac{\binom{m}{2}}{m!}\right)=
-\frac{m(m-1)}{(m!)^4}
\end{multline*}

Similarly, the coefficient of $n^{2m-1}$ in the
$\binom{n}{2m}$-term of $E(X_\tau^2)$ is
\[
-\frac{\binom{2m}{2}}{(2m)!}\frac{1}{(m!)^2}\binom{2m}{m}=
-\frac{m(2m-1)}{(m!)^4},
\]
so we only need to find the coefficient of the
$\binom{n}{2m-1}$-term of $E(X_\tau^2)$.

As we noted before, all subsets $S\subseteq[n]$ of the same size
(in our case, of size $2m-1$) are equivalent, so we may assume
$S=[2m-1]=\{1,2,\ldots,2m-1\}$. We want to find the number of
permutations $\rho$ of $S$ such that there exist subsets
$S_1,S_2\subseteq S$ of size $m$ for which we have $|S_1\cap
S_2|=1$ (so $S_1\cup S_2=S$) and $\rho(S)=\tau(S_1)\cup\tau(S_2)$.

Suppose that we want to choose $S_1$ and $S_2$ as above, together
with their positions in $S$, in such a way that the intersection
element $e$ is in the $i$th position in $\tau(S_1)$ and the $j$th
position in $\tau(S_2)$ (of course, $1\le i,j\le m$). Then $e$
occupies position $(i-1)+(j-1)+1=i+j-1$ in $S$. Hence, there are
$\binom{i-1+j-1}{i-1}$ ways to choose the positions for elements
of $\tau(S_1)$ and $\tau(S_2)$ on the left of $e$, and
$\binom{m-i+m-j}{m-j}$ ways to choose the positions for elements
of $\tau(S_1)$ and $\tau(S_2)$ on the right of $e$. On the other
hand, both $\tau(S_1)$ and $\tau(S_2)$ are naturally
order-isomorphic to $\tau$, hence, under that isomorphism $e$ maps
to $\tau(i)$ as an element of $S_1$ and to $\tau(j)$ as an element
of $S_2$. Since $e$ is the \emph{unique} intersection element, it
is easy to see that we must have
$e=(\tau(i)-1)+(\tau(j)-1)+1=\tau(i)+\tau(j)-1$ (exactly
$\tau(i)-1$ elements in $S_1$ and exactly $\tau(j)-1$ elements in
$S_2$, all distinct from those in $S_1$, must be less than $e$,
the rest of the elements of $S$ must be greater than $e$). There
are $\binom{\tau(i)-1+\tau(j)-1}{\tau(i)-1}$ ways to choose the
elements of $S_1$ and $S_2$ which are less than $e$, and
$\binom{m-\tau(i)+m-\tau(j)}{m-\tau(j)}$ ways to choose the
elements of $S_1$ and $S_2$ which are greater than $e$.

Thus, its positions in $\tau(S_1)$ and $\tau(S_2)$ uniquely
determine the position and value of the intersection element $e$;
there are $[i-1,j-1]_m$ ways to choose which other positions are
occupied by $\tau(S_1)$ and which ones by $\tau(S_2)$; and, there
are $[\tau(i)-1,\tau(j)-1]_m$ ways to choose which other values
are in $\tau(S_1)$ and which ones are in $\tau(S_2)$.

Now that we have chosen both positions and values of elements of
$S_1$ and $S_2$, we can produce a unique permutation $\rho(S)$ of
$S$ which satisfies our conditions above. Simply fill the
positions for $S_1$, resp. $S_2$, by elements of $\tau(S_1)$,
resp. $\tau(S_2)$, in the order in which they occur.

Since the total number of permutations of $S$ is $(2m-1)!$, the
coefficient of the $\binom{n}{2m-1}$-term of $E(X_\tau^2)$ is
\begin{multline*}
\frac{\sum_{i,j=1}^{m}{\binom{i-1+j-1}{i-1}\binom{m-i+m-j}{m-j}
\binom{\tau(i)-1+\tau(j)-1}{\tau(i)-1}\binom{m-\tau(i)+m-\tau(j)}{m-\tau(j)}}}
{(2m-1)!}=\\
=\frac{\sum_{i,j=1}^{m}{[i-1,j-1]_{m-1}[\tau(i)-1,\tau(j)-1]_{m-1}}}{(2m-1)!},
\end{multline*}
the coefficient of $n^{2m-1}$ in $V\!ar(X_\tau)$ is, by the
previous equations,
\begin{multline*}
\frac{\sum_{i,j=1}^{m}{[i-1,j-1]_{m-1}[\tau(i)-1,\tau(j)-1]_{m-1}}}{((2m-1)!)^2}-
\frac{m(2m-1)}{(m!)^4}+\frac{m(m-1)}{(m!)^4}=\\
=\frac{\sum_{i,j=1}^{m}{[i-1,j-1]_{m-1}[\tau(i)-1,\tau(j)-1]_{m-1}}}{((2m-1)!)^2}-
\frac{1}{(m!(m-1)!)^2}\ge0,
\end{multline*}
so we finally get
\[
\sum_{i,j=1}^{m}{[i-1,j-1]_{m-1}[\tau(i)-1,\tau(j)-1]_{m-1}}\ge
\left(\frac{(2m-1)!}{m!(m-1)!}\right)^2=\binom{2m-1}{m-1}^2,
\]
which is easily reducible to (\ref{eq-perm}) by $m\leftarrow m+1$,
then $\bar\tau(i)\leftarrow\tau(i+1)-1$. \qed

\bigskip

It seems, however, that a stronger form of our Theorem should be true,
namely, the following
\begin{conjecture} \label{conj1}
The \emph{strict} inequality holds in (\ref{eq-perm}) for all
$m>0$.
\end{conjecture}

This would imply that $V\!ar(X_\tau)$ has order $2m-1$ in $n$,
i.e. the standard deviation of $X_\tau$ is $1/2$ order smaller
than its expected value.

\begin{remark} \label{rem-covariance}
Similarly, the leading coefficient of the covariance
$Cov(X_{\tau_1},X_{\tau_2})$ is
\[
\frac{\sum_{i,j=1}^{m}{[i-1,j-1]_{m-1}[\tau_1(i)-1,\tau_2(j)-1]_{m-1}}}{((2m-1)!)^2}-
\frac{1}{(m!(m-1)!)^2},
\]
but
$\sum_{i,j=1}^{m}{[i-1,j-1]_{m-1}[\tau_1(i)-1,\tau_2(j)-1]_{m-1}}$
can be (and often is) less than $\binom{2m-1}{m-1}^2$.
\end{remark}

As of now, we only have some results on the sign of covariance for
small patterns. We hope to explore this topic further in
subsequent papers.

Note that the reversal map, $\tau(i)\mapsto\tau(m+1-i)$, the
complement map, $\tau(i)\mapsto m+1-\tau(i)$, preserve the
variance and covariance (we also make a note for the next section
that, for words $\tau\in[l]^m$, the reversal map is the same,
while the complement is $\tau(i)\mapsto l+1-\tau(i)$).

Considering symmetry classes of pairs of patterns (i.e.
equivalence classes with respect to reversal and complement), we
see that there are 8 classes of pairs of 3-letter patterns:
$\{123,123\}$, $\{132,132\}$, $\{123,132\}$, $\{132,213\}$,
$\{132,231\}$, $\{132,312\}$, $\{123,312\}$, $\{123,321\}$ (listed
in order of decreasing covariance). Of those, the first two pairs
obviously have a positive covariance, and of the remaining six,
only $\{123,132\}$ has a positive covariance.

Finally, denote the left-hand side and right-hand side of equation
(\ref{eq-perm}) by $L(m,\tau)$ and $R(m,\tau)$, respectively, and
let
\begin{gather*}
M^\ast(m)=\max_{\tau\in S_m}{(L(m,\tau)-R(m,\tau))},\\
M{\!}_\ast(m)=\min_{\tau\in S_m}{(L(m,\tau)-R(m,\tau))}.
\end{gather*}
It is not hard to see that $M^\ast(m)=L(m,id_m)-R(m,id_m)>0$,
where $id_m$ is the identity permutation of $\{0,1,\ldots,m\}$
(use Chebyshev's inequality, or dot product, or Cauchy-Schwarz
inequality). It would be interesting to characterize the
permutations $\tilde\tau_m$ such that
$M{\!}_\ast(m)=L(m,\tilde\tau_m)-R(m,\tilde\tau_m)$. We also make
the following conjecture.
\begin{conjecture} \label{conj2}
$\displaystyle{\exists\lim_{m\to\infty}{\frac{M{\!}_\ast(m)}{M^\ast(m)}}=0.}$
\end{conjecture}

\section{Patterns in words}

We now consider patterns contained in words, where repeated
letters are allowed both in the pattern and the ambient string.

\begin{theorem} \label{words}
Let $\tau$ be a map of $[0,m]=\{0,1,\ldots,m\}$ onto
$[0,l]=\{0,1,\ldots,l\}$. Then for any nonnegative integers $0\le
l\le m$ and any $\tau$ as above,
\begin{equation} \label{eq-words}
\sum_{i,j=0}^{m}{[i,j]_m [\tau(i),\tau(j)]_l} \ge
\frac{(2m+1)!(2l+1)!}{(m!)^2(l+1)!)^2}.
\end{equation}
\end{theorem}

\begin{remark} \label{rem-special}
Note that Theorem \ref{words} reduces to Theorem \ref{perm} when
$l=m$. Note also that, given $0\le l\le m$, Theorem \ref{words}
applies to $(l+1)!S(m+1,l+1)$ patterns $\tau$, where $S(m+1,l+1)$
is the Stirling number of the second kind.
\end{remark}

\begin{remark} \label{words-paths}
As in Theorem \ref{perm}, the left-hand side of Theorem
\ref{words} is the number of pairs $(P,Q)$ of northeast integer
lattice paths $P:(0,0)\to(i,j)\to(m,m)$ and
$Q:(0,0)\to(\tau(i),\tau(j))\to(l,l)$ over all $(i,j)\in[0,m]^2$.
\end{remark}

\noindent\textbf{Proof of Theorem \ref{words}}. The proof follows
the same outline as that of Theorem \ref{perm}, so we will use the
same notation as well. Again, it will be convenient to assume in
the proof that $\tau\in[l]^m$ is map of $[m]$ onto $[l]$ (i.e. use
$\{1,\dots,m\}$ instead of $\{0,1,\dots,m\}$ and $\{1,\dots,l\}$
instead of $\{0,1,\dots,l\}$) and, similarly, that the ambient
permutations $\sigma\in[k]^n$. Note that for any subset
$S\subseteq[n]$ of positions, the probability that the subsequence
of elements at positions in $S$, i.e. $\sigma(S)$, in a random
word $\sigma\in[k]^n$, is order-isomorphic to $\tau$ is
$\binom{k}{l}/k^m$. This is because $k^m$ is the total number of
subsequences of $m$ letters in $[k]$, $\tau$ has exactly $l$
distinct letters, and there are $\binom{k}{l}$ ways to choose $l$
distinct letters out of $k$. Hence, as in Theorem \ref{perm}, we
obtain
\[
E(X_\tau)=\frac{1}{k^m}\binom{k}{l}\binom{n}{m},
\]
which is a polynomial in $n$ and $k$. Therefore, the leading
coefficient of $E(X_\tau)$ as a polynomial in $n$ is
\[
[n^m]E(X_\tau)=\frac{1}{k^m}\binom{k}{l}\frac{1}{m!},
\]
so the leading coefficient of $E(X_\tau)^2$ is
\[
[n^{2m}]E(X_\tau)^2=\frac{1}{k^{2m}}\binom{k}{l}^2\frac{1}{(m!)^2}.
\]
However, as in the proof of Theorem \ref{perm}, we have that
$E(X_\tau^2)$ is a linear combination of $\binom{n}{2m-\ell}$,
$0\le\ell\le m$, with coefficients being polynomials in $k$ and
rational functions in $l,m$. A similar analysis shows that the
leading coefficient in $n$ of $E(X_\tau^2)$ is
\begin{multline*}
[n^{2m}]E(X_\tau^2)=\frac{1}{(2m)!}\left[\binom{n}{2m}\right]E(X_\tau^2)=\\
=\frac{1}{(2m)!}\binom{2m}{m}\binom{k}{l}^2\frac{1}{k^{2m}}=[n^{2m}]E(X_\tau)^2,
\end{multline*}
so $\deg_n(V\!ar(X_\tau))\le 2m-1$, and hence,
$[n^{2m-1}]V\!ar(X_\tau)\ge 0$.

As in the proof of Theorem \ref{perm}, we have that
\[
[n^{2m-1}]E(X_\tau)^2=2[n^{m-1}]E(X_\tau)[n^m]E(X_\tau)=
-\frac{m(m-1)}{(m!)^2}\binom{k}{l}^2\frac{1}{k^{2m}}
\]
and the coefficient of $n^{2m-1}$ in the $\binom{n}{2m}$-term of
$E(X_\tau^2)$ is
\[
-\binom{2m}{2}\frac{1}{(2m)!}\binom{2m}{m}\binom{k}{l}^2\frac{1}{k^{2m}}=
-\frac{m(2m-1)}{(m!)^2}\binom{k}{l}^2\frac{1}{k^{2m}}.
\]
The remaining summand in $[n^{2m-1}]V\!ar(X_\tau)$ is the
coefficient of $n^{2m-1}$ in the $\binom{n}{2m-1}$-term of
$E(X_\tau^2)$, i.e.
\begin{multline*}
\frac{1}{(2m-1)!}\left[\binom{n}{2m-1}\right]E(X_\tau^2)\\
-\frac{m(2m-1)}{(m!)^2}\binom{k}{l}^2\frac{1}{k^{2m}}
+\frac{m(m-1)}{(m!)^2}\binom{k}{l}^2\frac{1}{k^{2m}}\ge 0,
\end{multline*}
which is equivalent to
\[
\left[\binom{n}{2m-1}\right]E(X_\tau^2)\ge
\frac{(2m-1)!}{(m-1)!^2}\binom{k}{l}^2\frac{1}{k^{2m}}.
\]

As in the proof of Theorem \ref{perm}, it is easy to see that
$[\binom{n}{2m-1}]E(X_\tau^2)$ is equal to the probability that a
sequence $\rho\in[k]^{2m-1}$ is a union of two subsequences
order-isomorphic to $\tau$. Therefore, assume $[2m-1]=S_1\cup
S_2$, $\rho(S_1)\cong\tau\cong\rho(S_2)$. But then $S_1$ and $S_2$
have $m$ elements, so they intersect at a single element $e$.

Suppose that $e$ is at position $i$ in $S_1$ and at position $j$
in $S_2$. Then, as in the proof of Theorem \ref{perm}, there are
$\binom{i-1+j-1}{i-1}\binom{m-i+m-j}{m-i}=[i-1,j-1]_{m-1}$ ways to
choose which positions to the left and to the right of $e$ are in
$S_1$ and which ones are in $S_2$.

Suppose that $\rho$ contains $l+L$ distinct letters, then $0\le
L\le l-1$. Because of the positions of $e$ in $S_1$ and $S_2$, we
know that $e$ must map to $\tau(i)$ in $\rho(S_1)$ and to
$\tau(j)$ in $\rho(S_2)$ under our order-isomorphism. Suppose that
the value of $e$ in $\rho$ is $r$. Consider the $r-1$ letters in
$[l+L]$ which are less than $r$. Then
\[(r-1)-(\tau(j)-1)=r-\tau(j)\] of those occur only in $S_1$,
\[(r-1)-(\tau(i)-1)=r-\tau(i)\] occur only in $S_2$, and
\[(r-1)-(r-\tau(i))-(r-\tau(j))=\tau(i)+\tau(j)-1-r\] occur in both
$\rho(S_1)$ and $\rho(S_2)$. Similarly, of the $l+L-r$ letters in
$\rho$ which are greater than $r$,
\[(l+L-r)-(l-\tau(j))=L-r+\tau(j)\] occur only in $\rho(S_1)$,
\[(l+L-r)-(l-\tau(i))=L-r+\tau(i)\] occur only in $\rho(S_2)$, and
\[(l+L-r)-(L-r+\tau(i))-(L-r+\tau(j))=l-L+r-\tau(i)-\tau(j)\] occur
in both $\rho(S_1)$ and $\rho(S_2)$.

Thus, the number of sequences $\rho\in[k]^{2m-1}$ which are a
union of two subsequences order-isomorphic to $\tau$ is
\[
f(\tau,k)=\sum_{L=0}^{l-1}{\binom{k}{l+L}\sum_{r=0}^{l+L}{
\sum_{i,j=1}^{m}{[i-1,j-1]_{m-1}h(\tau,L,r,i,j)} } },
\]
where
\begin{multline*}
h(\tau,L,r,i,j)=\binom{r-1}{r-\tau(i),r-\tau(j),\tau(i)+\tau(j)-1-r}\times\\
\times\binom{l+L-r}{L-r+\tau(i),L-r+\tau(j),l-L+r-\tau(i)-\tau(j)}.
\end{multline*}

Hence, the probability that a sequence $\rho\in[k]^{2m-1}$ is a
union of two subsequences order-isomorphic to $\tau$ is
$f(\tau,k)/k^{2m-1}$, so we have
\[
\left[\binom{n}{2m-1}\right]E(X_\tau^2)=\frac{f(\tau,k)}{k^{2m-1}}\ge
\frac{(2m-1)!}{(m-1)!^2}\binom{k}{l}^2\frac{1}{k^{2m}},
\]
or, equivalently,
\[
kf(\tau,k)\ge\frac{(2m-1)!}{(m-1)!^2}\binom{k}{l}^2
\]
for all positive integers $k$ and all patterns $\tau\in[l]^m$. But
both sides of this inequality are polynomials in $k$ of degree
$2l$, hence the same inequality should hold for their leading
coefficients. The leading coefficient on the right is
\[
\frac{(2m-1)!}{(m-1)!^2}\frac{1}{(l!)^2}.
\]
On the left, $k^{2l}$ only occurs when $L=l-1$. But then
$\tau(i)+\tau(j)-1-r\ge 0$ and
$l-L+r-\tau(i)-\tau(j)=r+1-\tau(i)-\tau(j)\ge 0$, so
$r=\tau(i)+\tau(j)-1$, and hence
\begin{multline*}
h(\tau,L,r,i,j)=h(\tau,l-1,\tau(i)+\tau(j)-1,i,j)=\\
=\binom{\tau(i)+\tau(j)-2}{\tau(i)-1}\binom{2l-\tau(i)-\tau(j)}{l-\tau(i)}
=[\tau(i)-1,\tau(j)-1]_{l-1}.
\end{multline*}
Therefore,
\[
[k^{2l}](kf(\tau,k))=\frac{1}{(2l-1)!}\sum_{i,j=1}^{m}{[i-1,j-1]_{m-1}[\tau(i)-1,\tau(j)-1]_{l-1}},
\]
so
\[
\sum_{i,j=1}^{m}{[i-1,j-1]_{m-1}[\tau(i)-1,\tau(j)-1]_{l-1}}\ge\frac{(2m-1)!}{(m-1)!^2}\frac{(2l-1)!}{(l!)^2}.
\]
Now, letting $m\leftarrow m+1$, $l\leftarrow l+1$, then
$\bar\tau(i)\leftarrow\tau(i+1)-1$, we obtain the inequality
(\ref{eq-words}). \qed

\bigskip

Note that, for $l=0$ (which includes the case $m=0$), the
inequality (\ref{eq-words}) becomes an equality. We conjecture,
however, that the strict inequality holds if $l>0$, i.e. if $\tau$
is not a constant string.

As in the case of patterns in permutations, it would be
interesting to characterize the patterns $\tau\in[l]^m$, where the
difference between the two sides of (\ref{eq-words}) is minimal.

We also note that the covariance $Cov(X_{\tau_1},X_{\tau_2})$ of
patterns $\tau_1,\tau_2\in [l]^m$ is positive (resp. negative) if
\[
\sum_{i,j=1}^{m}{[i-1,j-1]_{m-1}[\tau_1(i)-1,\tau_2(j)-1]_{l-1}}-\frac{(2m-1)!}{(m-1)!^2}\frac{(2l-1)!}{(l!)^2}
\]
is positive (resp. negative). Hence, it would be interesting to
characterize pairs of patterns $\tau_1,\tau_2\in[l]^m$ based on
the sign of the covariance $Cov(X_{\tau_1},X_{\tau_2})$.

\medskip

\begin{center}
\textbf{Acknowledgements}
\end{center}

I am grateful to Herbert S. Wilf and Donald E. Knuth for their
helpful suggestions.


\begin{thebibliography}{9}
\bibitem {AAHHS} M.H. Albert, M.D. Atkinson, C.C. Handley, D.A. Holton, W.
Stromquist, On packing densities of permutations,
\textit{Electronic J. of Combinatorics} \textbf{9} (2002), \#R5.
\bibitem{AE} T. Amdeberhan, S.B. Ekhad, A condensed condensation proof of
a determinant evaluation conjectured by Greg Kuperberg and Jim Propp,
\textit{J. of Comb. Theory, Ser. A} \textbf{78} (1997), 169--170.
\bibitem{HLP} G.H. Hardy, J.E. Littlewood, G. P\'olya,
``Inequalities'', Cambridge University Press, Cambridge, 1934.
\bibitem{P} A. Price, Packing densities of layered patterns, Ph.D.
thesis, University of Pennsylvania, Philadelphia, PA, 1997.
\end{thebibliography}
\end{document}